\documentclass[12pt]{iopart}
\usepackage{amssymb,amsthm}
\usepackage{epsfig}
\usepackage{setstack}

\theoremstyle{plain}
\newtheorem{theorem}{Theorem}[section]

\theoremstyle{definition}
\newtheorem*{remark}{Remark}
\newtheorem{lemma}[theorem]{Lemma}

\input epsf

\begin{document}

\jl{8}

\title{Lower bounds for the energy in a crumpled elastic sheet -- A
minimal ridge}

\author{Shankar C Venkataramani \dag}

\address{\dag\ Department of Mathematics,
  University of Chicago, 5734 University Ave., Chicago, IL 60637,
  USA}

\ead{shankar@math.uchicago.edu}

\begin{abstract}
 We study the linearized F\"{o}pl -- von Karman theory of a long, thin
 rectangular elastic membrane that is bent through an angle $2
 \alpha$. We prove rigorous bounds for the minimum energy of this
 configuration in terms of the plate thickness $\sigma$ and the
 bending angle. We show that the minimum energy scales as
 $\sigma^{5/3} \alpha^{7/3}$. This scaling is in sharp contrast with
 previously obtained results for the linearized theory of thin sheets
 with isotropic compression boundary conditions, where the energy
 scales as $\sigma$.  
 \end{abstract}


\submitto{\NL}


\section{Introduction} \label{sec:intro}

Everyday experience tells us that thin elastic sheets {\em crumple},
when confined into a small volume, {\em e.g.} a sheet of paper
``confined'' by ones hands. Crumpling also plays an important role in
the mechanical behavior of packaging material and in the dissipation
of the energy of collisions by the ``crumple zones'' of
automobiles. Crumpling is therefore a problem of much intrinsic
interest, but understanding this behavior is complicated by the
complex morphology of a typical crumpled sheet.

Despite the complicated appearance of a crumpled sheet, the crumpling
phenomenon is in itself very robust. It is easily observed in thin
sheets made from a variety of materials, suggesting that it can be
studied using simplified or idealized models that capture the
essential features of thin elastic sheets. This approach leads one to
consider a crumpled sheet as a minimum energy configuration for a
simple elastic energy functional for thin sheets, {\em viz.} the
F\"{o}pl -- von Karman (FvK) energy. Using this approach of elastic
energy minimization, the crumpling response is now understood as a
result of the elastic energy of the sheet concentrating on a small
subset of the entire sheet \cite{pomeau,science.paper,eric}. The
energy in a crumpled sheet is concentrated on a network of thin
line-like creases (ridges) that meet in point-like vertices. Recent
work has resulted in quantitative understanding of both the vertices
\cite{maha,dcone_exp,vertex,BPCB00,MB02} and the ridges
\cite{lobkovsky,LobWit}. Scaling laws governing the behavior of
crumpled sheets have been obtained in the physics literature
\cite{pomeau,science.paper,lobkovsky,vertex}.
  
Minimum energy configurations for the FvK energy have also been
studied in the context of the {\em blistering problem}, {\em viz.} the
buckling of membranes as a result of isotropic compression along the
boundary.  The blistering problem is relevant to the delamination of
thin films that are chemically deposited at high temperatures, as well
as the mechanical behavior of micro-fabricated thin-film diaphragms
\cite{GiOrtz,GDeSOC02}.

There is a considerable body of mathematical work focused on the
blistering problem \cite{GiOrtz,KJ,Audoly,JS,JS1,BCDM,BCDM02}. Upper
and lower bounds have been obtained for approximations to the elastic
energy \cite{GiOrtz,KJ,JS1}, for the FvK energy \cite{JS,BCDM} and for
full three dimensional nonlinear elasticity \cite{BCDM02}. The FvK
energy and full three dimensional nonlinear elasticity give the same
scaling for the upper and the lower bounds. This yields a rigorous
scaling law for the energy of a blister. As we discuss below, the
scaling laws for the blistering problem are different from those for
the crumpling problem, even for the {\em same energy functional}. This
indicates that the energy minimizing configurations for the FvK energy
show an interesting dependence on the boundary conditions.

The minimization of the FvK elastic energy is an example of a
non-convex variational problem that is regularized by a singular
perturbation \cite{variational.review,sternberg}. It is well known
that this can lead to a variety of multiple-scale behaviors including
energy concentration and/or small scale oscillations in the minimizers
\cite{evans}.  

Multiple scale behaviors, both {\em microstructure} and {\em
singularities}, are ubiquitous in condensed matter systems. The
crumpling phenomenon appears to be a particularly simple and tractable
instance of multiple scale behavior. For this reason, there has been
much recent interest in physics literature about the nature of the
crumpling phenomenon \cite{nelson, Seung.Nelson, pomeau,
science.paper, origami}.  Asymptotic analysis \cite{lobkovsky} and
scaling arguments \cite{science.paper} show that the ridge energy
scales as
$$
{\cal E}_{r} \sim \sigma^{5/3} L^{1/3},
$$ where $\sigma$ is the thickness of the sheet, and $L$ is the length
of the ridge. Ben-Amar \etal \cite{pomeau} and Mahadevan \etal
\cite{maha} show that the energy of a $d$-cone scales as
$$
{\cal E}_{c} \sim \sigma^{2} \log(L/a),
$$ where $a$ is the radius of the core associated with the vertex. It
is clear that the ridge energy is asymptotically larger than the
vertex energy as $\sigma \rightarrow 0$. However, which of these two
energies is important for a given sheet depends on the relation
between the nondimensional thickness $\epsilon = \sigma/L$ of the
sheet, and the ``crossover thickness'' $ \epsilon ^* = \sigma/L$ which
is determined by setting ${\cal E}_r ={\cal E}_c$. $\epsilon^*$ in
turn depends on the values of the multiplicative constants for these
scaling laws. 

These multiplicative constants cannot be determined by scaling
arguments. Mahadevan \etal \cite{vertex} estimate the constant for the
vertex energy numerically using an ansatz for the shape of a sheet
near the vertex, and they find that
$$
{\cal E}_c \approx 100 \phi^2 \sigma^{2} \log(L/a),
$$ 
where $\phi$ is the complement of the tip angle of the cone
\cite{vertex,BPCB00}. Based on Lobkovsky's work \cite{lobkovsky},
Boudaoud \etal \cite{BPCB00} estimate the value of the constant in the
ridge energy as $1$. This implies, for the experiments in
Refs.~\cite{vertex}~and~\cite{BPCB00}, the vertex energy dominates the
ridge energy.

It would be useful to prove these scaling laws, and determine the
multiplicative constants rigorously, that is, in an ansatz-free
manner. In this paper, we propose a model problem that yields
structures analogous to a single ridge in a crumpled sheet. We prove a
rigorous lower bound (with a numerical value for the multiplicative
constant) for the elastic energy in our model problem. This is a step
toward rigorously proving the scaling law of Lobkovsky \etal for the
ridge energy \cite{science.paper}. We also discuss how the techniques
in this paper can be extended to prove similar results for a ``real''
crumpled sheet, as opposed to our model problem.

This paper is organized as follows -- In Sec.~\ref{sec:setup}, we
describe the problem of interest, set up the relevant energy
functional and determine the appropriate boundary conditions. We also
rescale the various quantities to a form that is suitable for further
analysis, and recast the problem in terms of the rescaled
quantities. In Sec.~\ref{sec:l_bnd}, we prove our main result, {\em
viz.} a lower bound for the elastic energy for our boundary
conditions. We present a concluding discussion in
Sec.~\ref{sec:discussion}.

\section{The Elastic energy} \label{sec:setup}

\begin{figure}[htbp]
  \begin{center}
\centerline{\epsfig{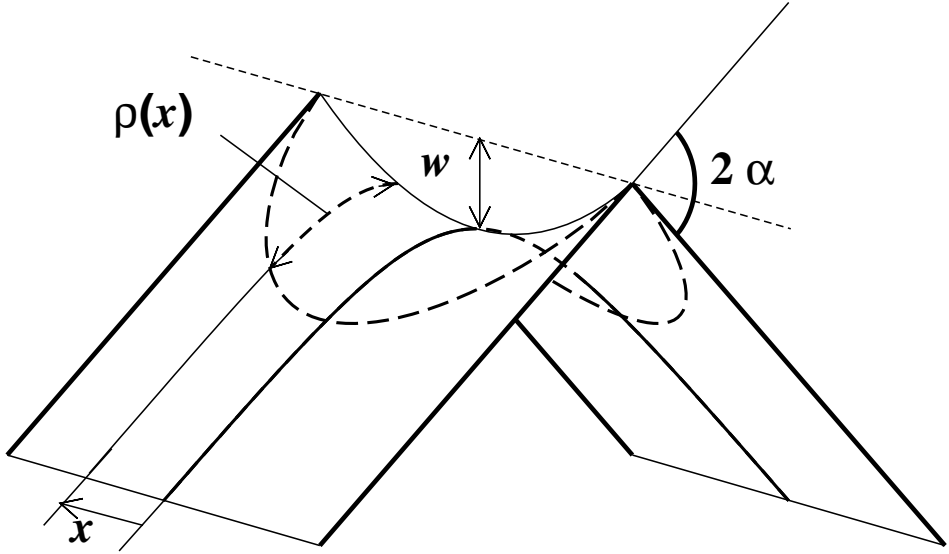}}
\end{center}
\caption{A minimal
    ridge. The boundary conditions are given by a frame (the thick
    solid lines) bent through an angle. The sheet is essentially flat
    outside the region bounded by the two dashed curves, and the bulk
    of the energy is concentrated in this region.}
\label{fig:ridge}
\end{figure}

We are interested in a {\em minimal ridge}, {\em i.e.}, the single
crease that is formed when a long rectangular elastic strip is bent
through an angle by clamping the lateral boundaries to a bent
frame. This situation is depicted in Figure~\ref{fig:ridge}.

\begin{figure}[htbp]
  \begin{center}
\centerline{\epsfig{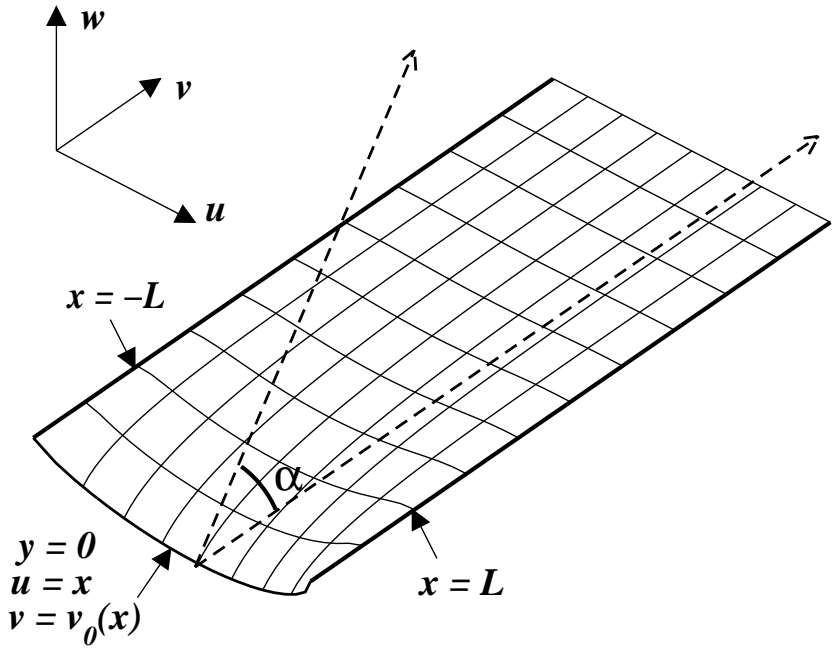}}
\end{center}
\vspace{0cm}
\caption{A schematic representation of our coordinate system and the
boundary conditions imposed on the sheet. The grid is given by the
lines $x = \mbox{constant}$ and $y = \mbox{constant}$. The two dashed
lines are the tangent to the curve $x = 0$ and the line $w = w_0(0), u
= 0$ respectively, and the angle between theses lines is $\alpha$.}
\label{fig:coords}
\end{figure}

From the symmetry of the problem, it is clear that we only need to
consider one half of the strip. This is represented schematically in
Figure~\ref{fig:coords}. We will use (the material) coordinates
$(x,y)$ on the reference half strip $|x| \leq L$, $0 \leq y \leq
\infty$. Also, we associate a $(u,v,w)$ coordinate systems in space,
so that as $y \rightarrow \infty$, the half-strip is asymptotically in
the $w = 0$ plane, as depicted in Figure~\ref{fig:coords}. The grid in
the figure is generated by the lines $x = \mbox{constant}$ and $y =
\mbox{constant}$ that are straight in the reference (material)
coordinates. $w$ represents the out of plane deformation, and the
in-plane coordinate are chosen so that the $u$ and the $v$ axes are
asymptotically parallel to the $x$ and $y$ axes respectively as $y
\rightarrow \infty$. Since the sheet is bent through an angle $2
\alpha$, as $y \rightarrow - \infty$, the sheet will lie in the plane
$w = v \, \tan(2 \alpha)$. The symmetry of the two halves implies that
the the line $y = 0$ maps into the plane $w \, \tan \alpha + v = 0$,
which bisects the angle between the planes $w = 0$ (the asymptote as
$y \rightarrow \infty$) and $w = v \, \tan(2 \alpha)$ (the asymptote
as $y \rightarrow -\infty$). The symmetry between the two halves also
implies that the tangent to the lines $x = \mbox{constant}$ at $y = 0$
should be perpendicular to the plane $w \, \tan \alpha + v =
0$. Consequently $u_y = 0$ and $w_y = v_y \, \tan \alpha$ at $y = 0$.

A mathematically justified way to obtain the elastic energy of the
deformed sheet is to treat the sheet as a three dimensional (albeit
thin) object and use a full nonlinear three dimensional elastic energy
functional for the energy density. This approach however does not take
advantage of the ``thinness'' of the sheet. In particular, we would
like to treat the thin sheet as a two dimensional object. This will
greatly reduce the complexity of the problem. The derivation of
reduced dimensional descriptions of thin sheets has a long
history. There is a classical theory for thin elastic sheets built on
the work of Euler, Cauchy, Kirchoff, F\"{o}pl and Von Karman
\cite{love,landau_elastic}.

Many modern authors have considered the problem of deriving a reduced
dimensional theory \cite{ciarlet2} from three dimensional elasticity
in a mathematically rigorous fashion through $\Gamma$ -- convergence
\cite{DeGDM,DalMaso}. $\Gamma$ -- convergence is the appropriate
notion for convergence in variational problems. In the context of thin
sheets, roughly speaking, finding the $\Gamma$ -- limit amounts to
identifying an appropriate two dimensional energy functional whose
minimizers give the limiting behavior of the minimizers of the full
three dimensional energy in the limit the thickness $\sigma
\rightarrow 0$. This problem hasn't been solved in its entirety. 

There are two situations for which reduced dimensional theories have
been derived rigorously as $\Gamma$ -- limits of full three
dimensional elasticity. Membrane theories \cite{LeDR93,LeDR95,LeDR96}
are applicable in situations where the stretching is essentially
uniform through the thickness of the sheet and bending theories
\cite{FJM} are applicable in situations where the the sheet is
essentially unstretched, and all the elastic energy is due to strains
that are first order in the thickness of the sheet.

Neither of these theories are applicable for the minimal ridge. For
the minimal ridge, both the stretching (membrane) energy and the
bending energy are important. In fact, a ridge is a result of the
competition between these two energies. So we turn to the classical
F\"{o}pl -- von Karman {\em ansatz} for an appropriate definition of
the elastic energy \cite{love,landau_elastic}. Although the FvK energy
is not rigorously derived from three dimensional elasticity, it can be
thought of as the sum of the membrane and the bending energies that
have been derived rigorously in different limiting situations. We will
further assume that, for the minimal ridge, the deflections
$|w|,|u-x|$ and $|v-y|$ are small compared to the natural length scale
$L$. The strains are of the order of the square of the maximum
deflection divided by $L$ and they are small. After some rescaling,
the energy of the deformed sheet is given by the linearized FvK energy
\begin{eqnarray}
 {\mathcal I} & = & \int\left[(u_x + w_x^2 - 1)^2 + \frac{1}{2} (v_x +
u_y + 2 w_x w_y)^2 + (v_y + w_y^2 - 1)^2\right] dx dy \nonumber \\ & &
+ \sigma^2 \int (w_{xx}^2 + 2 w_{xy}^2 + w_{yy}^2 ) \, dx dy.
\label{eq:unscaled} 
\end{eqnarray}
where $x$ and $y$ are reference coordinates on the sheet, $u$ and $v$
are in-plane coordinates, $w$ is the out of plane displacement and
$\sigma$ is the thickness of the sheet. The integrand for the first
integral is given by the squares of the linearized strains,
$$ 
\gamma_{xx} = u_x + w_x^2 - 1, \quad \gamma_{xy} = \gamma_{yx} = w_x w_y
+ \frac{1}{2}(v_x + u_y), \quad \gamma_{yy} = v_y + w_y^2 - 1.
$$ The blistering of thin films is also described by the elastic
energy in (\ref{eq:unscaled}). A similar energy also describes
multiple scale buckling in {\em free} elastic sheets ({\em i.e.}
sheets that are not forced through the boundary conditions) that are
not intrinsically flat \cite{eran}.

The difference between the blistering problem and a minimal ridge is
in the boundary conditions, which we describe below. Since the strains
are assumed to be small, $v_y + w_y^2 \approx 1$. If the bending
half-angle $\alpha
\ll 1$, so that $\tan \alpha \approx
\alpha$, we get $w_y = \alpha \mbox{ at } y = 0$.  Since the
deformation goes to zero far away from the bend, we have
$$
|u - x| \rightarrow 0, |v - y| \rightarrow 0, w \rightarrow 0 \quad
 \mbox{ as } y \rightarrow \infty,
$$
Also, the lateral boundaries at $x = \pm L$ are clamped to the
frame. Therefore,
$$
u = x, v = w = 0, \mbox{ at } x = \pm L.
$$ Note that the appropriate boundary condition for the minimal ridge
at $y = 0$ is a free boundary condition, subject to the symmetry
requirement $w + v \tan \alpha = 0$. We are going to replace this free
boundary condition with a Dirichlet boundary condition
$$
v = v_0(x), w = w_0(x) \quad \mbox{ at } y = 0,
$$
where we will leave the functions $v_0(x)$ and $w_0(x)$ unspecified,
except for a size condition. Defining 
$$
a = \max_{[-L,L]}(|v_0(x)|,|w_0(x)|),
$$ 
we impose the size condition by requiring that $a$ be ``small''.

\subsection{Rescalings} 

The relevant small parameter in the problem is the non-dimensional
thickness of the sheet, $\epsilon = \sigma/L$. Following Lobkovsky
\cite{lobkovsky}, we introduce the rescaled coordinates and
displacements by
$$
x = L X, \quad \quad y = \sigma^{1/3} L^{2/3} Y,
$$
and 
$$
w = \sigma^{1/3} L^{2/3} W, \quad \quad v = y + \sigma^{1/3}
L^{2/3} V, \quad \quad u = x + \sigma^{2/3} L^{1/3} U.
$$
Since $\sigma, L,x, y, u, v,w$ all have dimensions of a length, it
is clear that the rescaled quantities $X,Y,U,V,W$ are all
dimensionless. With these rescalings, the dimensionless energy $I =
\sigma^{-5/3} L^{-1/3} {\mathcal I}$ is given by
\begin{eqnarray}
\fl I(U,V,W) & = & \int \left[(U_X + W_X^2)^2 +
\frac{\epsilon^{-2/3}}{2} (V_X + U_Y + 2 W_X W_Y)^2 +
\epsilon^{-4/3}(V_Y + W_Y^2)^2 \right] \nonumber
\\ 
& & + \left[ W_{YY}^2 + 2 \epsilon^{2/3} W_{XY}^2 +\epsilon^{4/3}
W_{XX}^2 \right]dX dY.
\label{eq:scaled} 
\end{eqnarray}
Our quest for rigorous scaling results for the energy ${\mathcal I}$
reduces to the following -- Show that the rescaled energy $I$, of a
minimizer $(U^*,V^*,W^*)$, is bounded above and below by positive 
constants {\em uniform} in the dimensionless thickness parameter
$\epsilon$, as $\epsilon \rightarrow 0$.

Setting 
$$ W_0(X) = \sigma^{-1/3} L^{-2/3} w_0(LX), \quad V_0(X) =
\sigma^{-1/3} L^{-2/3} v_0(LX),
$$
the rescaled quantities satisfy the boundary conditions
$$
V = V_0(X), W = W_0(X) , W_Y = \alpha \quad \mbox{ at } Y = 0,
$$
and 
$$
U \rightarrow 0, V \rightarrow 0, W \rightarrow 0 \quad \mbox{ as } Y
 \rightarrow \infty,
$$
We have the lateral boundary conditions
$$
U = V = W = 0 \quad \mbox{ at } X = \pm 1.
$$ 
Also, the deformation at $Y = 0$ satisfies $|V_0(X)| \leq A$ and
$|W_0(X)| \leq A$ where $ A = \sigma^{-1/3} L^{-2/3} a$.

\section{Lower Bound} \label{sec:l_bnd}

In this section, we prove a lower bound for the linearized Elastic
energy in Eq.~(\ref{eq:unscaled}), provided that the length scale $a$
associated with the boundary conditions satisfies a size condition.

\begin{theorem} ${\mathcal I}(u,v,w)$ is as defined in Eq.~(\ref{eq:unscaled}). 
There is a constant $b > 0$ such that, for all $u \in H^1, v \in H^1$
and $w \in H^2 \cap W^{1,4}$ satisfying the boundary conditions
\begin{eqnarray}
|u - x| \rightarrow 0, |v - y| \rightarrow 0, w \rightarrow 0 \quad
 \mbox{ as } y \rightarrow \infty, \nonumber \\
u = x, v = w = 0, \mbox{ at } x = \pm L, \nonumber \\
v = v_0(x), w = w_0(x) \quad \mbox{ at } y = 0, \nonumber
\end{eqnarray}
and the size condition
$$
\max_{[-L,L]}(|v_0(x)|,|w_0(x)|) = a \leq b \sigma^{1/3} L^{2/3} \alpha^{2/3},
$$
we have the lower bound
$$
{\mathcal I}(u,v,w) \geq \frac{2}{5} \alpha^{7/3} \sigma^{5/3} L^{1/3}.
$$
\end{theorem}

\begin{remark} We will prove the theorem by proving the scaled version of 
the statement, {\em viz.}, the rescaled boundary conditions along with
the rescaled size condition
$$
\max(|V_0(X)|,|W_0(X)|) = A \leq b \alpha^{2/3},
$$ 
imply that
$$
I \geq \frac{2 \alpha^{7/3}}{5}.
$$
\end{remark}
\begin{remark} The hypothesis for the theorem includes a size condition on 
the displacement at $y = 0$. This is somewhat unsatisfying, as it is
not {\em a priori} obvious that a ridge in a ``real'' crumpled sheet
will satisfy this condition.
\end{remark}

In our search for a lower bound, we can assume $I(U,V,W) < \infty$
w.l.o.g. From Eq.~(\ref{eq:scaled}), it follows that for $\sigma > 0$,
$W \in H^2$. By the standard trace theorems, the boundary conditions
$W(X,0) = W_0(X)$ and $W_Y(X,0) =
\alpha$ are therefore assumed pointwise for almost every $X$.

Since $U = 0$ at $X = \pm 1$, it follows that $ \int_{-1}^1 (U_X +
W_X^2) dX = \int_{-1}^1 W_X^2 dX$. Using this together with Jensen's
inequality we get
\begin{eqnarray}
I(U,V,W) & \geq & \int \left[(U_X + W_X^2)^2 + W_{YY}^2 \right]dX dY
\nonumber \\ & \geq & \int_0^{\infty} \left[ \frac{1}{2} \left(
\int_{-1}^1 W_X^2 dX \right)^2 + \int_{-1}^1 W_{YY}^2 dX \right] dY 
\label{eq:onlyw}
\end{eqnarray}
Thus the functional
$$
E(W) = \int_0^{\infty} \left[ \frac{1}{2} \left(\int_{-1}^1 W_X^2 dX 
\right)^2 + \int_{-1}^1 W_{YY}^2 dX \right] dY 
$$ 
bounds $I(U,V,W)$ from below. Hence it suffices to obtain a lower
bound for $E$ with the given hypothesis to prove the theorem.

Let $E_b$ and $E_s$ denote the quantities
$$
E_b = \int W_{YY}^2 dX dY,
$$
and
$$
E_s = \int_0^{\infty}\frac{1}{2} \left(\int_{-1}^1 W_X^2 dX \right)^2
dY.
$$ Although $E_b$ and $E_s$ are only lower bounds for the ``true''
bending and the stretching energies $I_b$ and $I_s$, we will
nonetheless call $E_b$ and $E_s$ the bending and the stretching energy
for convenience.

For every $X$, we define
$$
\rho(X) = \left[\int_0^{\infty} W_{YY}^2(X,Y) dY\right]^{-1},
$$
and for every $Y$, we define
$$
\tau(Y) = \int_{-1}^1 W_X^2 dX.
$$ Since $E < \infty$, $\rho(X) > 0$ (a.e.) and $\tau(Y) < \infty$
(a.e.). 

For any function $f$ depending on $X$ and $Y$, let $\|f\|_Y^2$ denote
$\int_{-1}^{1} |f(X,Y)|^2 dX$, so that,
$$
\tau(Y) = \|W_X\|^2_Y.
$$ $\rho(X)$ is a ``local'' (in $X$) measure of the bending energy,
and $[\rho(X)]^{-1}$ can be thought of as the bending energy density
in $X$ that is obtained by integrating out the $Y$ dependence. We can
also think of $\rho(X)$ as the length scale associated with the
bending energy as a function of $X$, {\em viz.}, we expect that the
bending energy density in $Y$ decays rapidly for $Y/\rho(X) \gg 1$
(See Fig.~\ref{fig:ridge}). Likewise, $[\tau(Y)]^2$ is a local (in
$Y$) measure of the stretching energy. In terms of $\rho(X)$ and
$\tau(Y)$, the bending an the stretching energies are given by
$$
E_b = \int_{-1}^{1} \frac{1}{\rho(X)} dX \quad \quad E_s =
\int_0^{\infty}\frac{1}{2} [\tau(Y)]^2 dY.
$$

We now begin our proof of the theorem. The idea behind the proof is to
show that the stretching energy $E_s$ can be bounded from below by a
negative power of the bending energy $E_b$, so that the total energy
$E_s + E_b$ tends to $+\infty$ as $E_b \rightarrow 0$ and $E_b
\rightarrow \infty$. This ensures the existence of a positive lower
bound for $E$ (and consequently also for $I$).

\begin{lemma} For every $Y$, we have the inequality
$$
\|W\|_Y^2 \geq \alpha^2 Y^2 - 2 Y^3 E_b  - 4 A^2. 
$$
\end{lemma}
\begin{proof}
If $\rho(X) > 0$, $W(X,.)$ is a $C^1$ function by the Sobolev Embedding
theorem and the boundary conditions imply that $W(X,0) = W_0(X)$ and
$W_Y(X,0) = \alpha$. Consequently, for such an $X$,
$$
W(X,Y) = W_0(X) + \alpha Y + \int_0^Y W_{YY}(X,\xi) (Y - \xi) d \xi.
$$
We will estimate $W^2(X,Y)$ from this equation using the elementary
inequalities
$$
|a - b|^2 \geq (1 - \delta) |a|^2 - \frac{1 - \delta}{\delta} |b|^2,
$$
and 
$$
|a + b|^2 \leq (1 + \delta) |a|^2 + \frac{1 + \delta}{\delta} |b|^2,
$$
for all $\delta > 0$. 
By our hypothesis on the boundary conditions,
$$
|W_0(X)|^2 \leq A^2.
$$
By Jensen's inequality
\begin{eqnarray}
\left(\int_0^Y W_{YY}(X,\xi) (Y - \xi) d \xi\right)^2 & \leq & Y 
\int_0^Y W^2_{YY}(X,\xi) (Y - \xi)^2 d \xi \nonumber \\
& \leq & Y^3 \int_0^Y W^2_{YY}(X,\xi) d \xi \nonumber \\
& \leq & \frac{Y^3}{\rho(X)} \nonumber
\end{eqnarray}
Combining these estimates, we get
$$
W^2(X,Y) \geq (1 - \delta_1)\alpha^2 Y^2 - \frac{1 - \delta_1}{\delta_1}
\left[(1 + \delta_2) \frac{Y^3}{\rho(X)}  + \frac{1 + \delta_2}{\delta_2} 
A^2 \right],
$$ for all positive $\delta_1$ and $\delta_2$. In particular, setting
$\delta_1 = 1/2$ and $\delta_2 = 1$ yields
\begin{equation}
W^2(X,Y) \geq \frac{1}{2}\alpha^2 Y^2 - \frac{2 Y^3}{\rho(X)} - 2 A^2.
\label{eq:local_bound} 
\end{equation} 
Integrating this inequality in $X$ we obtain
\begin{equation}
\|W\|_Y^2 \geq \alpha^2 Y^2 - 2 Y^3 E_b  - 4 A^2, 
\label{eq:integral_bound}
\end{equation}
proving the lemma.
\end{proof}

Our proof is based on demonstrating that a small bending energy
$E_b$ will lead to a large stretching energy. This idea is quantified
by the following lemma where we use Eq.~(\ref{eq:integral_bound}) to
estimate the stretching energy $E_s$ from below in terms of the
bending energy $E_b$.

\begin{lemma} Let
$$
\mu = \left(\frac{4 E_b A}{\alpha^3}\right)^2.
$$
There is a constant $\mu^* > 0$ such that, if $\mu < \mu^*$, the
stretching energy $E_s$ satisfies
$$
E_s \geq \frac{\alpha^{14}}{5 \cdot 729 \cdot E_b^5}.
$$
\label{lem:essential}
\end{lemma}

\begin{proof}
We have,
$$
\tau(Y) = \|W_X\|_Y^2 \geq \frac{\pi^2}{4} \|W\|_Y^2,
$$
by the Poincare Inequality. Eq.~(\ref{eq:integral_bound}) now yields
$$
\tau(Y) \geq \frac{\pi^2}{4}\left[\alpha^2 Y^2 - 2 Y^3 E_b - 4 A^2
\right].
$$ We will now extract the appropriate scaling dependence of $\tau(Y)$
and $E_s$ on $\alpha$ and $E_b$. Setting $\alpha^2 Y^2 = 2 Y^3 E_b$,
we deduce that a characteristic scale $\tilde{Y}$ for $Y$ is given by
$$
\tilde{Y} = \frac{\alpha^2}{2 E_b}.
$$
Rescaling $Y$ in terms of $\tilde{Y}$, we obtain
$$
\tau(Y) \geq \frac{\alpha^6 \pi^2}{16
E_b^2}\left[\left(\frac{Y}{\tilde{Y}}\right)^2 \left(1 -
\frac{Y}{\tilde{Y}}\right) - \mu \right],  
$$ 
where $\mu$ is as defined above, {\em i.e.}
$$
\mu = 
\left(\frac{4 E_b A}{\alpha^3}\right)^2.
$$ Making the change of variables $Y = z \tilde{Y}$, and using
$\tau(Y) \geq 0$, we see that
$$
E_s \geq \int_0^{\infty} \frac{1}{2} [\tau(Y)]^2 dY \geq
\frac{\pi^4 \alpha^{14} K(\mu)}{512 E_b^5},
$$
where 
$$
K(\mu) = \int_0^{1} \left[ \max(z^2(1-z) - \mu,0)\right]^2 dz.
$$
$K(\mu)$ is clearly a continuous function of $\mu$ and $K(0) =
1/105$. Since $\pi^4 K(0)/512 > 1/(5 \cdot 729)$, there is an $\mu^* >
0$ such that for all $\mu < \mu^*$, $\pi^4 K(\mu)/512 > 1/(5 \cdot
729)$. The lemma follows.
\end{proof}
We can now prove the theorem.
\begin{proof}
Let $b = \frac{5}{8} \sqrt{\mu^*}$. The hypothesis imply
$$
A < \frac{5}{8} \sqrt{\mu^*}
\alpha^{2/3}.
$$
If $E_b \geq 2 \alpha^{7/3}/5$, there is nothing to prove. Therefore,
we can assume that $E_b < 2 \alpha^{7/3}/5$. Then, it follows that
$$
\mu = 
\left(\frac{4 E_b A}{\alpha^3}\right)^2 < \mu^*.
$$
Consequently $E_s \geq \alpha^{14}/(5 \cdot 729 \cdot
E_b^5)$. Minimizing $E_s + E_b$, we see that
$$
I \geq E_s + E_b \geq \frac{\alpha^{14} }{5 \cdot 729 \cdot E_b^5} + E_b \geq
\frac{2 \alpha^{7/3}}{5}.
$$ 
The theorem follows by ``undoing'' our rescaling to express the
results in terms of $u,v,w,a$ and ${\mathcal I}$.
\end{proof}

\begin{remark}
In the appendix, we show that 
$$
K(\mu) \geq \frac{1}{105} - \frac{\mu}{6}.
$$ Therefore, we can choose $\mu^* = 0.048$ and $b = 0.13$. We see
that $b$ is not exceedingly small. Rather it is $O(1)$. Also, it is
not the best constant for this theorem, and we can get a better value
by optimizing our choices for the constants ($\delta_1$ and $\delta_2$
which we set to be $1/2$ and $1$ respectively).
\end{remark}

\section{Discussion} \label{sec:discussion}

We have proved a rigorous lower bound for the elastic energy of a
ridge that is consistent with the results obtained by Lobkovsky \etal
\cite{science.paper} using scaling arguments, and by Lobkovsky
\cite{lobkovsky} using matched asymptotics. In order to prove these
scaling laws rigorously, we also need analogous upper bounds that are
consistent with the same scaling. This approach has been used
successfully for a variety of other variational problems
\cite{KM,JS,BCDM}. The upper bound is usually obtained by constructing
a test solution that yields the bound. One is often guided by scaling
arguments in the construction of the appropriate test solution. This
is in contrast to the lower bounds, where one needs to obtain
functional analysis type inequalities that captures the competition
between distinct energies in the problem ({\em e.g.}
lemma~\ref{lem:essential}). It is this competition that determines the
scaling behavior of the problem.

We have not constructed the upper bounds, because we believe that they
will follow from the scaling ansatz in Lobkovsky's work
\cite{lobkovsky}. We also believe that the upper bound will scale in
the same manner as the lower bound, thereby giving us a rigorous
scaling law for the energy of a single ridge.

As we remark before proving the theorem, our result is not directly
applicable to the minimal ridge since we have the extra hypothesis 
$$
a \leq b \sigma^{2/3}\alpha^{7/3} L^{1/3},
$$
We need to replace this restriction with a free boundary condition at
$y = 0$ subject to the symmetry requirement $w \, \tan \alpha + v = 0$,
and the boundary condition $w_y = v_y \, \tan \alpha$. 

Despite this caveat, we claim that the analysis in this paper captures
the essential features of the minimal ridge problem, namely the
scaling in lemma~\ref{lem:essential}, and consequently the scaling law
for the energy of the ridge. We will show that this is indeed the case
by investigating the problem with the ``true'' boundary conditions
elsewhere.

A harder problem is to show that the scaling laws also hold for a real
crumpled sheet, where the forcing is not through clamping the
boundaries to a frame, but rather through the confinement in a small
volume. In this case, there are interesting global geometric and
topological considerations, some of which are explored in
Refs.~\cite{immersion_thm} and \cite{high_d}. As Lobkovsky and Witten
\cite{LobWit} argue, the boundary condition that the deformation goes
to zero far away from the ridge implies that the ridges do not
interact with each other significantly. The ridges can be considered
the {\em elementary excitations} of a crumpled sheet. In our quest for
upper and lower bounds, this translates to the statement that the test
solutions for the upper bound can be constructed by piecing together
local solutions near each ridge. Thus the upper bound is obtained
relatively easily. The hard part is to show that the confinement
actually leads to the formation of ridges, and that the competition
between the bending and the stretching energy for this situation has
the same form as in lemma~\ref{lem:essential}. In this context, we
expect that global topological considerations, as well as the non
self-intersection of the sheet will play a key role in the analysis,
as they do in the analysis of elastic rods (one dimensional objects)
\cite{GMSvdM}.

As we note above, the blistering problem is described by the same
elastic energy (Eq.~(\ref{eq:unscaled})), but with different boundary
conditions. Our results show an interesting contrast with results for
the blistering problem. Ben Belgacem {\em et al.} have shown that
\cite{BCDM}, for an isotropically compressed thin film, the energy of
the minimizer satisfies
$$
c \lambda^{3/2} \sigma L \leq {\mathcal I}
\leq C \lambda^{3/2} \sigma L,
$$
where $L$ is a typical length scale of the domain, and $\lambda$ is
the compression factor. A construction for the upper bound strongly
suggests that the minimizers develop an infinitely branched network
with oscillations on increasingly finer scales as $\sigma
\rightarrow 0$. In contrast, our results indicate that the energy of a 
minimal ridge satisfies
$$
c \sigma^{5/3} L^{1/3} \leq {\mathcal I}
\leq C \sigma^{5/3} L^{1/3},
$$ and the energy concentrates in a region of width
$\sigma^{1/3}L^{2/3}$. This shows that the nature of the solution of
the variational problem for the elastic energy in (\ref{eq:unscaled})
depends very strongly on the boundary conditions. In particular the
very nature of the energy minimizers is different for the two problems
-- For the blistering problems, as $\sigma \rightarrow 0$ the
minimizers develop a branched network of folds refining towards the
boundary. The minimizers therefore display the problem of small scale
{\em oscillations} \cite{tartar,evans}. The minimal ridge problem on
the other hand displays the {\em concentration} phenomenon
\cite{tartar,evans} as $\sigma \rightarrow 0$, with the energy
concentrating on a region of width $\sim \sigma^{1/3} L^{2/3}$. It
would be interesting to explore this issue further, and in particular,
to understand the mechanisms by which the boundary conditions
determine the nature of the minimizer.

\appendix

\section*{Appendix}

\setcounter{section}{1}

In this appendix, we prove the inequality
$$
K(\mu) = \int_0^{1} \left[ \max(z^2(1-z) - \mu,0)\right]^2 dz \geq
\frac{1}{105} - \frac{\mu}{6}.  
$$
\begin{proof}
Let $f \geq 0$. Then, for $f \geq \mu$, we have
$$
\left[ \max(f - \mu,0)\right]^2 = f^2 - 2 \mu f + \mu^2 \geq f^2 - 2 \mu f.
$$
For $0 \leq f \leq \mu$, we have
$$
\left[ \max(f - \mu,0)\right]^2 = 0 \geq f^2 - 2 \mu f.
$$
Using these inequalities in the definition of $K(\mu)$ we obtain
$$
K(\mu) \geq \int_0^{1} \left[z^2(1-z)\right]^2 - 2 \mu z^2 (1 - z) dz
\geq \frac{1}{105} - \frac{\mu}{6}.  
$$
\end{proof}

\ack

I wish to thank Brian DiDonna, Bob Kohn, L. Mahadevan, Stefan
M\"{u}ller and Tom Witten for helpful and enlightening
conversations. This work was started during a visit to the Max-Planck
Institute for Mathematics in the Sciences in Leipzig, and I would like
to thank the MPI for their hospitality. This work was supported in
part by the National Science Foundation through its MRSEC program
under Award Number DMR-9808595, and by a NSF CAREER Award
DMS-0135078. This work is also supported in part by a Research
Fellowship from the Alfred P. Sloan Jr. Foundation.

\end{document}